\documentclass{article}

\usepackage[english]{babel}

\usepackage[letterpaper,top=2cm,bottom=2cm,left=3cm,right=3cm,marginparwidth=1.75cm]{geometry}

\usepackage{amsmath}
\usepackage{graphicx}
\usepackage[colorlinks=true, allcolors=blue]{hyperref}

\usepackage{amssymb}

\usepackage{amsthm}
\newtheorem{proposition}{Proposition}

\title{An improved lower bound to Erdos' problem concerning products of distances for fixed diameter}
\author{Nat Sothanaphan\footnote{natsothanaphan@gmail.com}}
\date{December 16, 2025}

\begin{document}
\maketitle

\begin{abstract}
Erdos, Herzog and Piranian asked whether, for $n$ points in the plane with fixed diameter (maximum distance between points), an arrangement of a regular $n$-gon maximizes their product of all pairs of distances. Recently, it was discovered that, for every even $n \geq 4$, a regular $n$-gon is not a maximizer. However, the discovered improvement turns out to be very small. Indeed, for a fixed diameter of $2$, let $\Delta$ be the square of the product of all pairs of distances (the ``square'' is here due to connections with polynomial discriminants). Then, for a regular $n$-gon, $\Delta = n^n$ for even $n$. The discovered arrangements have proven $\Delta = (1+o(1))n^n$ thus far, and it was not known whether one can have $\Delta \geq C n^n$ for some $C > 1$ and all sufficiently large even $n$. In this note, we show that indeed $\liminf_{n\to\infty} \Delta_{\max}/n^n > 1.037$ for even $n$ which settles this conjecture. Other arrangements with higher conjectured $\Delta/n^n$ values are in fact known, but we have not been able to obtain proofs that they have large products of distances. Finally, no arrangements such that $\Delta/n^n \to \infty$ are known and we do not know whether they exist.
\end{abstract}

\section{Introduction}

Erdos, Herzog and Piranian \cite{ErdosHerzogPiranian1958} asked the following question. Let $n \geq 2$ and $f$ be a monic polynomial with roots $z_1,\dots,z_n \in \mathbb{C}$. Let the absolute value of the discriminant of $f$ be
$$ \Delta(f) := \prod_{i \neq j} |z_i - z_j| = \prod_{i < j} |z_i - z_j|^2, $$
that is, the square of product of distances between all pairs of roots. Assume further that the set of roots has diameter $2$:
$$ \max_{i,j} |z_i - z_j| = 2. $$
Then, it is asked whether $\Delta(f)$ is maximized when the $z_i$'s form a regular $n$-gon. For even $n$, this configuration gives $\Delta = n^n$. For odd $n$, this configuration gives $\Delta \sim e^{\pi^2/8} n^n$ (the extra factor is from the fact that the configuration is slightly bigger than the unit circle in order for the largest diagonal to have length $2$), where $\sim$ means the ratio of both sides tends to $1$ as $n \to \infty$.

In this note, we prove the following result. Let $\Delta_{\max}$ denote the largest value of $\Delta(f)$ according to the ``fixed diameter of $2$'' condition. Then:

\begin{proposition}\label{prop:lowerBd}
We have, for even $n$,
$$ \liminf_{n \to \infty} \frac{\Delta_{\max}}{n^n} \geq C \approx 1.0378 > 1, $$
where $C$ is the constant given by
$$ C := \exp \left(-I \frac{\pi^2}{128}\right), $$
and $I$ is the definite integral
$$ I := \int_0^{2\pi} \int_0^{2\pi} \frac{((1-\frac{x}{\pi})e^{ix/2} - (1-\frac{y}{\pi})e^{iy/2})^2 (e^{ix} + e^{iy})}{(e^{ix} - e^{iy})^2} \,dx\,dy \approx -0.481436 < 0. $$
In fact, there is an explicit construction for each even $n$ with
$$ \Delta = C n^n (1 + O(1/n)). $$
\end{proposition}

The ``strange'' form of this constant follows from our proof technique utilizing vector fields.

We note that it is likely this constant $C$ is not best possible, but such results appear difficult to prove, and this one happened to be the first proven case. We invite future investigators to improve upon this lower bound.

Finally, it is not currently known whether $\Delta_{\max} / n^n \to \infty$ as $n \to \infty$ for even $n$.

The author thanks Stijn Cambie, ChatGPT, and Quanyu Tang for discussion (in alphabetical order).

\section{Previous literature}

As far as the author is aware, there is limited literature on this question. In \cite{ErdosHerzogPiranian1958} itself, it was proved that if the set $\{ z: |f(z)| \leq 1 \}$ is connected, then $\Delta(f) \leq n^n$. This motivated the question itself. Pommerenke \cite{Pommerenke1961} worked on this question and obtained the upper bound $\Delta_{\max} \leq 2^{4(n-1)} n^n$. Here the factor $\Delta_{\max}/n^n$ is exponential; and we don't know of any better upper bound at the moment. So there is a large gap between the regular $n$-gon values of $n^n$ for even $n$ or $e^{\pi^2/8}n^n(1+o(1))$ for odd $n$ and the known upper bound.

Recently, there is renewed progress on the Erdos problems website, where the problem is numbered at \#1045 \cite{BloomErdosProblem1045}. Here's a little summary. Cambie, Hu and Tang showed that, for even $n$, if you move one diameter segment radially a little, the diameter of $2$ is preserved and the product of distances increases a little. But the improvement is modest: around $\Delta = (1 + O(1/n^2))n^n$.

An improvement here was obtained by Cambie, Dong and Tang, where they considered a regular Reuleaux $n/2$-gon (curve of constant width) along with arc midpoints. This improves the lower bound to $\Delta = (1 + O(1/n))n^n$ for $n$ even. Still relatively modest! It was then asked whether $\Delta_{\max} = (1 + o(1))n^n$ for $n$ even, expecting this to be false. But so far no proof has been obtained.

Cambie, Dong and Tang then obtained a more sophisticated construction for $n=6k$ which consists of $6$ circular arcs and obeys some symmetry (it is $D_3$ symmetry). Here, numerics suggests $\Delta/n^n$ tends to a value around $1.30$, which would clear the $(1+o(1))n^n$ bar. However, this configuration resists closed-form expressions and they (along with the author) have so far been unable to show the improved lower bound of $C n^n$ for some $C>1$. Hence, that conjecture remained unresolved.

In the present note, the author finally resolves this conjecture by a different construction.

Finally, the case of $n$ odd remains without any improvement; no construction better than the regular $n$-gon has so far been discovered in that case.

\section{The construction with new lower bound}

The author extended the idea of ``moving diameters radially'', as has been discussed on the Erdos problems website \cite{BloomErdosProblem1045}.

\subsection{Description of the construction}

We note that the Reuleaux-based configuration described in the previous section is in fact an instance of this moving-diameter idea. There, we alternate between ``pushing'' and ``pulling''.

However, it is actually better to continuously deform the ``push amounts'' rather than alternating them like that. This led to the following construction. Throughout, $n$ must be even.

Let:

$$ \delta_k = c \left( \frac{1}{n} - \frac{4k}{n^2} \right), z_k = (1 + \delta_k) e^{2\pi i k/n}, z_{k + n/2} = -(1 - \delta_k) e^{2\pi i k/n}, \quad k = 0,1,\dots,\frac{n}{2}-1. $$

Here $c = c_{\max}(n)$ is the maximum quantity that maintains the diameter of $2$. Numerically, $c_{\max}(n) \to \pi^2/4$ as $n \to \infty$. This will also be proven later. The author obtained numerical evidence that $\Delta / n^n \to 1.037... > 1$ for this construction, but was initially unable to show it just like the case of $n=6k$ construction of Cambie, Dong and Tang (where the value is instead around $1.30$ for that other case).

The idea here is that we're moving the diameter at angle $0$ outward by $c/n$. To be consistent, the diameter at angle $\pi$ must be moved by $-c/n$. Then we linearly interpolate the diameters in between them.

\subsection{Delicateness of the construction}

I want to draw attention to one aspect of such constructions that may seem surprising and is certainly non-obvious. See how I linearly interpolate the push amounts. It may be thought that the manner in which we interpolate doesn't matter. This is in fact not the case, and shows why the proof has been difficult!

If instead I change the formula in the previous subsection to
$$ \delta_k = \frac{c}{n} \cos\left( \frac{2 \pi k}{n} \right),$$
then actually numerics suggests $\Delta / n^n \to 1$ very rapidly. Thus, this would not defy the lower bound of $\Delta = (1+o(1))n^n$ as we might hope.

Currently, there is no good understanding on the part of the author on how to determine which interpolations work.

Similarly, the mentioned $n=6k$ construction with circular arcs is conceivably ``sensitive'' in similar ways; this would explain why the proof of large $\Delta$ has been difficult to obtain in that case as well.

\section{Proof of new lower bound}

We now present proof of Proposition \ref{prop:lowerBd}.

\subsection{Vector field}

The idea is that instead of looking at the configuration as a fixed entity, we imagine it ``flowing'' according to some vector field $v$.

It is important that this vector field $v$ is independent of $n$, the number of points. This happens to work out by coincidence for my construction (I did not intend that there is a ``fixed'' expression here; it happened during the proof). The vector field in question turns out to be:
$$ z' = v(z) := \left( 1 - \frac{2}{\pi} |\arg z| \right) \frac{z}{|z|}, \quad \arg z \in [-\pi, \pi). $$
Note this is not defined for $z=0$, but it's ok because we only need it near the unit circle. Also, the vector field is continuous (due to $|\arg z|$), but it's not differentiable. It is however Lipschitz continuous near the unit circle; this fact will turn out to be crucial.

The construction turns out to be exactly the following. We begin with the $n$th roots of unity (forming a regular $n$-gon). Then move these points by the above differential equation from $t=0$ to $t=O(1/n)$. (The exact ending $t$ is $c_{\max}/n$, which will turn out to be $O(1/n)$. This will be proven later in a more precise form.)

\subsection{Taylor expansion}

Now I define $v_j := v(z_j)$. Also define the ``differentials'':
$$ \rho_{ij} := \frac{v_i - v_j}{z_i - z_j}. $$
Due to Lipschitz continuity mentioned in the previous subsection, $\rho_{ij}$ is bounded (independent of $n$). This is very important!

Now we notice
$$ \log \Delta = \sum_{i,j} \log |(z_i - z_j) + (v_i - v_j)t| = \text{const} + \sum_{i,j} \log |1 + \rho_{ij} t|. $$
We do Taylor expansion in the following way.
$$ \log |1 + \rho t| = \text{Re} \log (1 + \rho t) = \text{Re} \left( \rho t - \frac{\rho^2 t^2}{2} + \frac{\rho^3 t^3}{3} - \frac{\rho^4 t^4}{4} + \dots \right). $$
There is a symmetry $v(-z) = v(z)$, which means $z_{j+n/2} = -z_j$ while $v_{j+n/2} = v_j$. This implies that corresponding to $\rho_{ij}$ there will also be a term that is its negative. So in the summation above in $\log \Delta$, the odd-power terms in the Taylor expansion will cancel; we're only left with even-power terms.

Extend the notation of $\rho_{ij}$ to also be able to be evaluated at any $t$; here $z_j(t)$ is the $z_j$ that has flowed for time $t$ and again $v_j(t) = v(z_j(t))$ (the notation for $v_j$ is for completeness only because $v_j$ is constant in the flow in our configuration anyway). Recall the Lagrange form of the remainder in Taylor's theorem. This means:
$$ \log \Delta = \log (n^n) - \frac{1}{2} \sum_{i,j} \text{Re}(\rho_{ij}^2) t_{\max}^2 - \frac{1}{4} \sum_{i,j} \text{Re}(\rho_{ij}(t_{\text{intermediate}})^4) t_{\max}^4 $$
where $t_{\max}$ is the ending time and $t_{\text{intermediate}}$ is some value in $[0, t_{\max}]$.

This form will suffice for our estimation!

\subsection{Derivation of bound}

We continue using the formula for $\log \Delta$ derived in the previous subsection.

\subsubsection{The second derivative term}

First, let's tackle the main $t_{\max}^2$ term. We must sum up all of $\rho_{ij}^2$. Using symmetry, I simplified it to
$$ \sum_{i,j} \rho_{ij}^2 = 4 \sum_{0 \leq i,j < n/2} \frac{(v_i - v_j)^2(z_i^2 + z_j^2)}{(z_i^2 - z_j^2)^2} = 4 \cdot \frac{n^2}{4} \left( \frac{I}{4\pi^2} + O\left(\frac{1}{n}\right) \right), $$
where $I$ is the integral in Proposition \ref{prop:lowerBd}. To get the integral $I$ in that form, note that the $z_j^2$ are the $n/2$th roots of unity, and then use some algebra to transform our finite sum into a Riemann sum approximating the integral $I$. The expression $I/(4\pi^2)$ there is the average value of the integrand.

Now I must explain how I get that $O(1/n)$ error term. First, due to the important fact that $\rho_{ij}$ is bounded, the integrand in $I$ and each term in our finite sum is bounded independent of $n$. So the first observation is that there are only $n/2(n/2-1)$ terms in the summation rather than $n^2/4$; but the missing terms only have combined size at most $O(n)$, so it fits into the error term there. From now on, assume we truly have $n^2/4$ sampling points arranged as $n/2 \times n/2$ grid. The integrand of $I$ is in fact smooth (the $(x-y)^2$ factor in the denominator also appears in the numerator), so the approximation on this grid has error $O(1/n)$ per term. All of this come out to be as claimed.

It also turns out that $I$ is a real number; this is a consequence of another of our symmetry that $v(\bar{z}) = \overline{v(z)}$. The fact that it is a \emph{negative} real number, however, is central to making $C > 1$ (looking back to the statement of Proposition \ref{prop:lowerBd}). But this negativity cannot be derived as an easy consequence of anything. Instead, Wolframalpha evaluated $I$ to be a negative real number. This is related to why this method is able to prove large $\Delta$ at all: it identifies a critical integral $I$ that just happens to satisfy our desired inequality. Again, the author does not currently have good understanding on which vector fields $v$ will produce the corresponding $I$ that comes out negative as we needed, only that this particular construction works as intended.

\subsubsection{The fourth derivative term}

This is the $t_{\max}^4$ term. Note that, with $t = t_{\text{intermediate}}$,
$$ |\sum_{i,j} \rho_{ij}(t)^4| = \left| \sum_{i,j} \frac{\rho_{ij}^4}{(1 + \rho_{ij} t)^4} \right| \leq \sum_{i,j} \frac{|\rho_{ij}|^4}{(1 - |\rho_{ij} t|)^4} = O(n^2). $$
This is because of the crucial fact that $\rho_{ij}$ is bounded and we have $t=O(1/n)$ so that $|\rho_{ij} t| = o(1)$.

This is actually a smaller error than what is needed; in the main formula of $\log \Delta$ this makes the final fourth-derivative contribution only $O(1/n^2)$ (when the error combined from all sources should be $O(1/n)$).

\subsubsection{The time to flow}

Finally we must estimate $t_{\max}$ itself. This is the maximum $t$ that maintains diameter of $2$. The correct rate here is
$$ t_{\max} = \frac{\pi^2}{4n} (1 + O(1/n)). $$
There are probably many ways to show this. We will present one way.

We solve the equation $|z_i(t) - z_j(t)| = 2$. This is $\log |z_i - z_j| + \log |1 + \rho_{ij} t| = \log 2$. Let $\rho = \rho_{ij}$. Expanding the $\log |1 + \rho t|$ to second derivative, we get
$$ t = - \log |(z_i - z_j)/2| (1 + O(t)) / \text{Re}(\rho). $$
Direct calculation gives $\text{Re}(\rho) = 1 - (\text{angle between }z_i,z_j)/\pi$. By some further calculations, we find that $t$ is constrained most strongly when $z_i,z_j$ are almost opposite to each other: the angle between them is $\pi - 2\pi/n$. Also we have $O(t) = O(1/n)$. Plugging all those in gives the desired rate for $t_{\max}$.

\subsubsection{Finishing}

Gathering the three error estimates, we arrive at
$$ \log \Delta = \log(n^n) - \frac{1}{2} \cdot \frac{n^2 I}{4 \pi^2} \left( 1 + O\left(\frac{1}{n}\right) \right) \cdot \frac{\pi^4}{16n^2} \left( 1 + O\left(\frac{1}{n}\right) \right) + O\left(\frac{1}{n^2}\right). $$
This gives the exact form in the statement of Proposition \ref{prop:lowerBd} and finishes the proof.

\bibliographystyle{alpha}
\bibliography{references}

\end{document}